\documentclass{amsart}
\usepackage{amsmath}
\usepackage{amsfonts}
\usepackage{amssymb}
\usepackage{epsfig}
\usepackage{amscd}
\usepackage[all]{xy}
\usepackage{color}
\usepackage{graphicx}

\begin{document}

\begin{center}

{\bf\LARGE Hyperbolic flowers}

\vspace{0.5cm}

{\sc\Large Maria Trnkova} 

\smallskip

{\it Department of Mathematics, University of California, Davis, CA 95616, USA.}
{\tt mtrnkova@math.ucdavis.edu}
\end{center}

\medskip

\begin{center}
{\it Abstract}
\end{center}

\vspace{-0.1cm}

\noindent Crochet models of a hyperbolic plane is a popular educational tool as they help to visualize complicated objets in hyperbolic geometry. We present another way how to make crochet models when we view them as a part of a triangulated hyperbolic plane. We  also provide a model of a cylinder in a hyperbolic space. This approach helps to understand various properties of hyperbolic geometry that are demonstrated in the paper: a sum of angles and a relation between edges and angles in a hyperbolic triangle, tiling of a hyperbolic plane, ratio of the circumference to the radius of a hyperbolic disc and even Nash-Kuiper embedding theorem. Oriented on students learning basics of Riemannian geometry. 


\section{Introduction}

Bill Thurston introduced several ways how to visualize hyperbolic plane and see some of its properties. He proposed several ways of making it with paper triangles and annuli, fabric pentagons \cite{Th}. Later he continued this work with Kelly Delp and they made more complicated surfaces \cite{Delp}.  This approach is described in details by Kathryn Mann in her notes ``Do-it-yourself: Hyperbolic Geometry" and contains many exercises of the topic  \cite{Mann}.
Another very esthetic models where introduced by Taimina \cite{HT, Taimina} - how to crochet hyperbolic disks. It attracted a lot of attention and it was followed by contributions of many other artists and even appeared in a TED talk \cite{TED}. It  is well described in a Craft magazine \cite{Craft}.

Our approach is inspired by these models but we make one step further and instead of crochet solid disks we make them triangulated.  The main difference is that these new models allow to explain some intrinsic properties of hyperbolic geometry. We also introduce a crochet model of a hyperbolic cylinder. In this note we want to use the advantages of triangulated models to demonstrate some well known results of hyperbolic geometry: (i) relation between edges and angles in a hyperbolic triangle, (ii) tiling of a hyperbolic plane, (iii) circumference growth of a disc, (iv) Nash-Kuiper embedding theorem, (v) visualization of a cylinder in hyperbolic space. The crochet scheme is given in the figure 9. The models can be used in classes of different levels: from high school to graduate courses to demonstrate special features of negatively curved spaces as well as the comparison with Euclidean spaces and surfaces. 

\section{Constructions with equilateral triangles}

Let us take a collection of equilateral triangles of the same size and connect them along their edges without shifting or stretching. We classify different cases based on the number of triangles meeting at each vertex. We can visualize the process by taking the actual triangles made out of wood, paper or plastic and try to glue them together. The simplest case is two triangles, and if we glue them along each edge we immediately get a {\it triangle pillow}. Following the same strategy for three triangles meeting at the same point we get a {\it tetrahedron}, for four triangles an {\it octahedron} and for five triangles a {\it icosahedron}. 

\begin{figure}[h]
\begin{center}
    \includegraphics[scale=0.17]{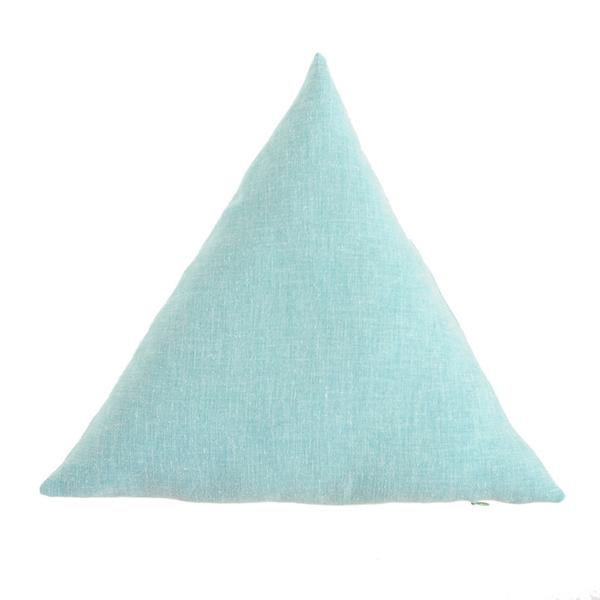}
    \includegraphics[scale=0.1]{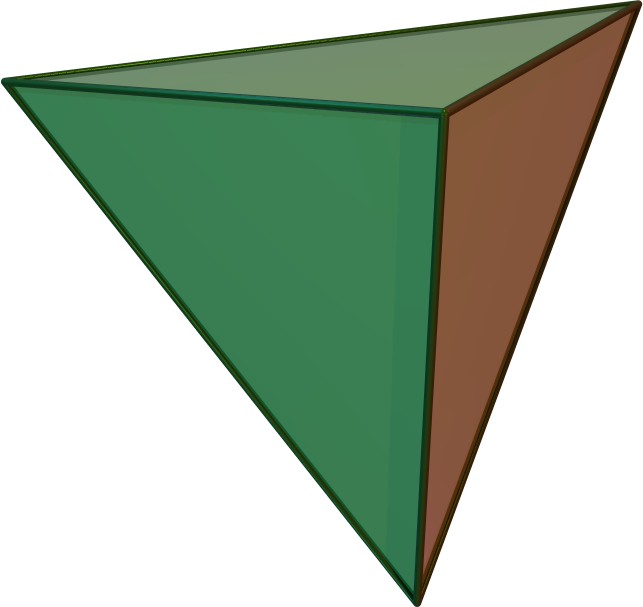}
  \includegraphics[scale=0.08]{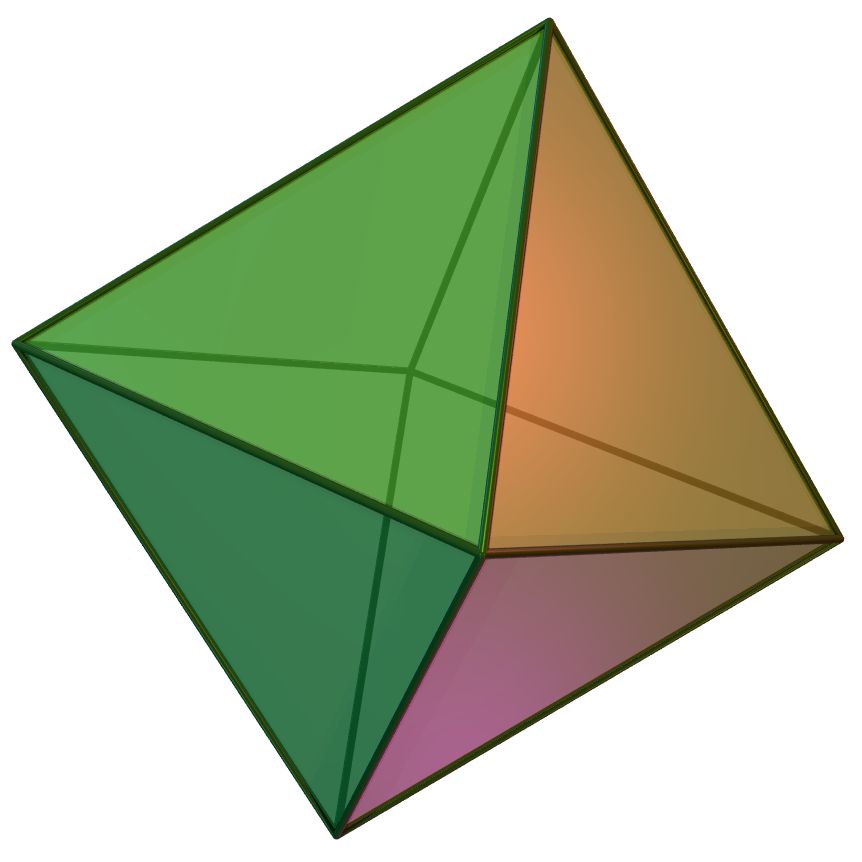}
  \includegraphics[scale=0.2]{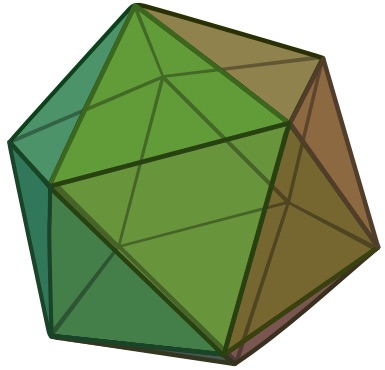}\\
  \caption{Gluing equilateral triangles with pattern of two, three, four and five at each vertex. Pictures courtesy of Wikipedia.}
\end{center}
\end{figure}

For the next case of six triangles we realize that all triangles which meet in one vertex must lie in a plane and form a hexagon. If we continue adding six triangles to the previously obtained vertices we tile the whole Euclidean plane.

\begin{figure}[h]
  \begin{center}
  \includegraphics[scale=0.65]{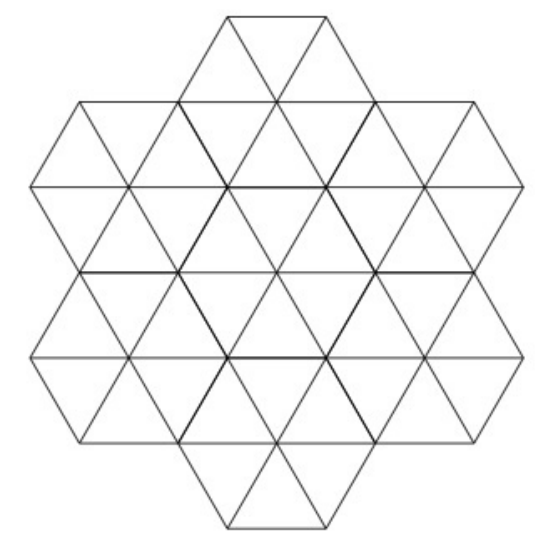}
  \end{center}
  \vspace{-0.2cm}
  \caption{Tiling of Euclidean plane with triangles}
\end{figure}

We saw that if the number of triangles was less than six we needed a finite number of triangles and the resulting object was a closed object of {\it positive curvature}. In the case of six triangles we got a full plane and the number of triangles was infinite. It has zero curvature. If we go even further and demand having seven triangles at each vertex something very interesting happens. We can start with a single vertex, add seven triangles and get a polygon which cannot be flatten. Nevertheless we continue to add more triangles to other vertices, and keep this process going. We quickly realize that we again need an infinite number of triangles to satisfy the condition of having seven triangles at each vertex. Unlike in the six triangles case the resulting object does not look like a flat plane -- it is a {\it hyperbolic plane} $\mathbb H^2$, and it has a negative curvature. For a better visualization we can choose some material and glue the triangles together. This would be harder to do with a firm material like paper, wood or plastic, but this can be nicely done with a yarn. We apply the idea above for creating pieces of hyperbolic plane and  crochet these models. The triangulation of the models allow to observe some geometric properties. The crochet models in figure 3 are examples of hyperbolic planes.

\section{Topology and Geometry}

\subsection{Tilings of a hyperbolic plane}
It is very well known that only few regular polygons can tile the Euclidean plane: triangles, squares and hexagons. Note that size of each of the shapes is not fixed. The situation is very different in the hyperbolic case. There are {\it infinitely} many ways how to tile a hyperbolic plane $\mathbb H^2$ with {\it regular} polygons. The most natural way is to tile it with equilateral triangles.  Unlike in the Euclidean case there is a freedom in a number of equilateral triangles at each vertex. Once we choose this number it determines the size of an equilateral triangle, it means that the tiling cannot be scaled.
If there are $k$ triangles around each vertex of a hyperbolic plane then every angle of an equilateral triangle must be $2\pi/k$. The smallest number $k$ of equilateral triangles around one vertex that gives a hyperbolic structure is seven. The interesting fact is that the sum of angles in each triangle is $6\pi/k$ which depends on $k$, $k\geq7$ and it is always less than $\pi$.  For Euclidean triangles the sum of internal angles is always $\pi$, and this difference is the main intrinsic property of a hyperbolic plane. Our two models (right two objects in figure 3) have seven equilateral triangles around each vertex. They approximate a disk area in $\mathbb H^2$ and also demonstrate that this pattern could be uniquely extended to a tiling of the plane. The crochet model looks the same around different vertices or triangles.

\begin{figure}
  \begin{center}
  \includegraphics[scale=0.4]{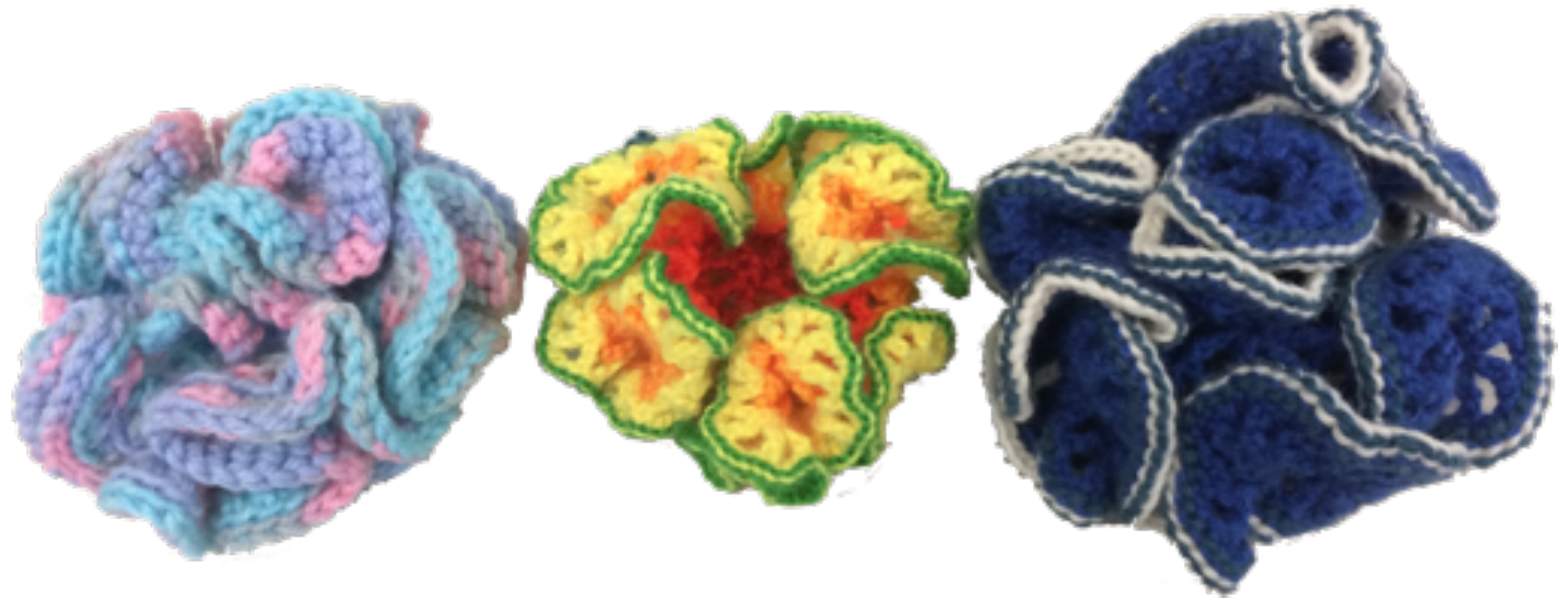}
  \end{center}
  \vspace{-0.2cm}
  \caption{The crochet model on the left is a solid hyperbolic disk. Two models on the right are  triangulated hyperbolic surfaces, approximating hyperbolic disks. They have contour made in a contrast color to visualize their boundaries and better hold shapes.}
  \label{flower3}
\end{figure}

There are different models of a hyperbolic plane  $\mathbb H^2$. We will consider the projective disc model  $\mathbb D^2$ also known as Beltrami-Klein model \cite{Th, S}. It is an open unit disc with the following metric $d_{D^2}(P,Q)$ in $D^2$,

$$
\cosh d_{D^2}(P,Q) = \frac{1-P\cdot Q}{\sqrt{1-||P||^2}\sqrt{1-||Q||^2}},
$$

where $P, Q$ are vectors in $D^2$, $||P||$ is a Euclidean norm of $P$. Geodesics in this model look like chords of the disc $D^2$. 

\begin{figure}[h]
  \begin{center}
  \includegraphics[scale=0.23]{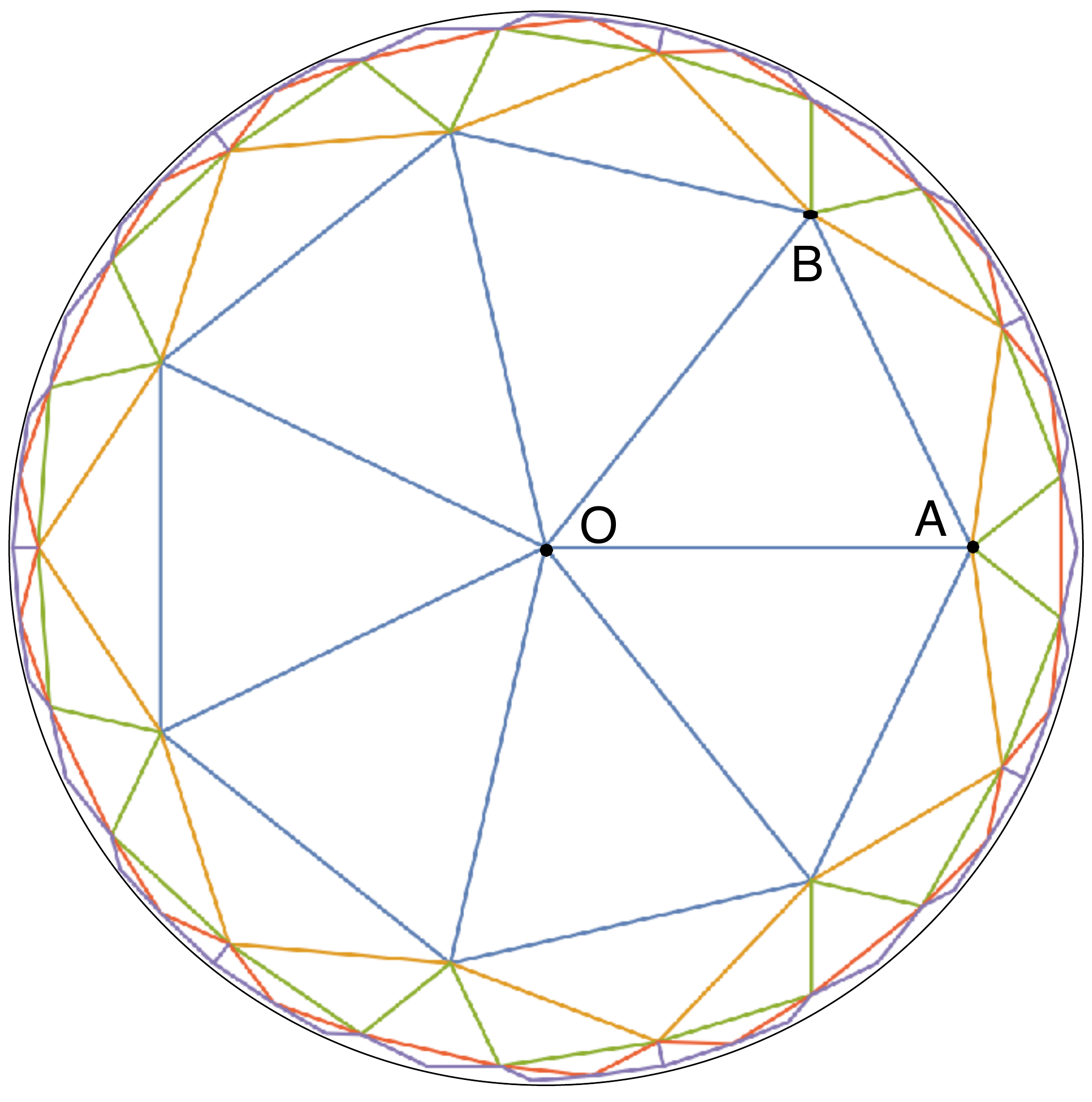} \qquad
    \includegraphics[scale=0.21]{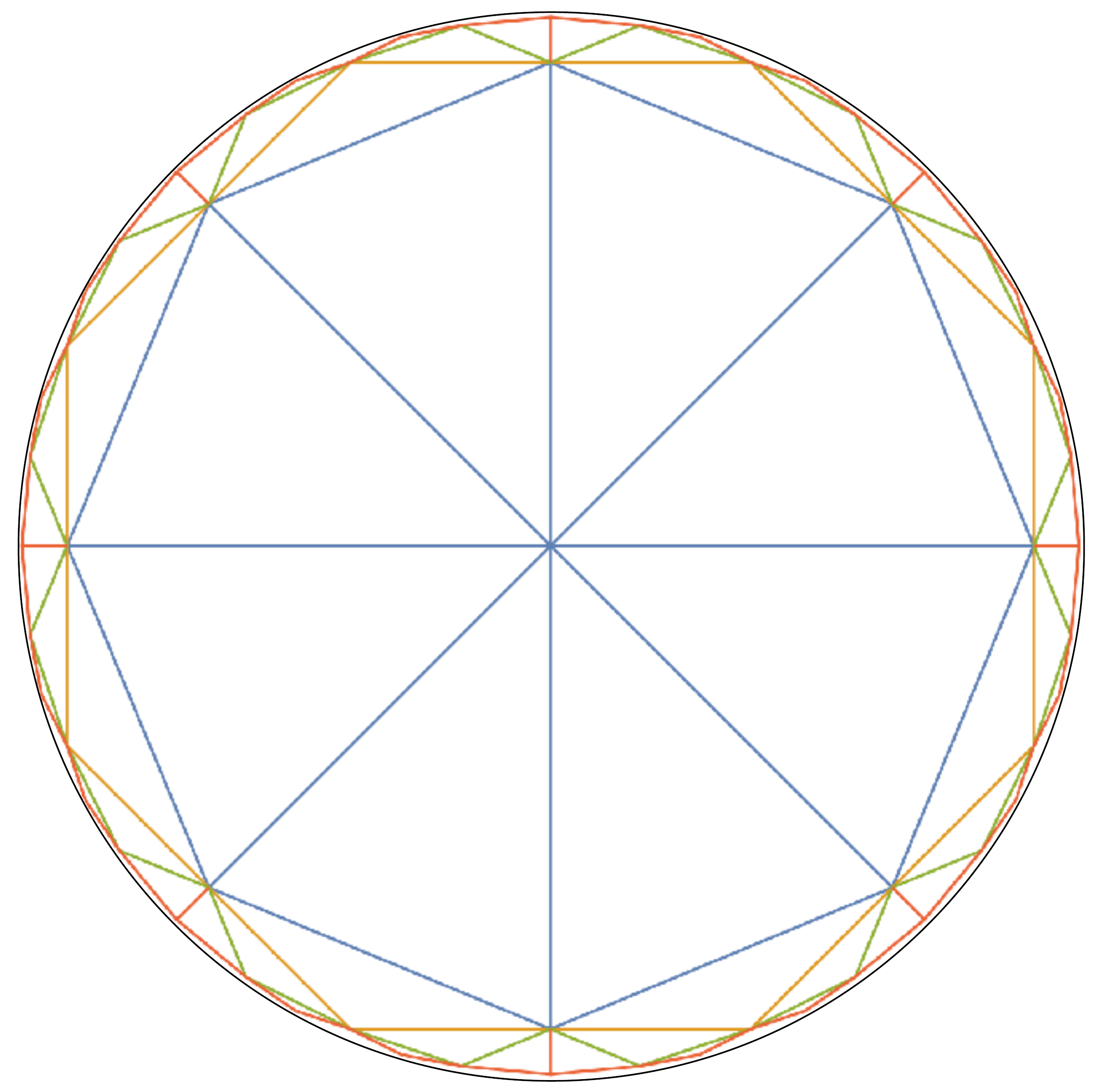}
    \end{center}
    \vspace{-0.2cm}
  \caption{Several equilateral hyperbolic triangles in the projective disc model that tile the hyperbolic plane with seven, respectively eight triangles at each vertex.}
  \label{fig1}
\end{figure}

We consider the case when seven equilateral triangles meet at each vertex. Every angle of such triangle is $2\pi/7\approx 51.4^\circ$ although it does not look so in the projective disc model in figure \ref{fig1}. We use the above distance formula to compute coordinates and edge length of the equilateral triangle. Let one vertex be at the origin $O(0, 0)$ and the other one $A$ on $x-$axis. Then  the $A$ and the third vertex $B$ should have coordinates $(a_1,0)$ and $(b_1, b_2)$ that satisfy the equation 
$$\cosh d_{D^2}(O,A)=\cosh d_{D^2}(O,B)=\cosh d_{D^2}(A,B),$$ where points $O, A, B$ can be viewed as vectors. Rewrite these equations once again

$$
\frac{1-O\cdot A}{\sqrt{1-||O||^2}\sqrt{1-||A||^2}}=\frac{1-O\cdot B}{\sqrt{1-||O||^2}\sqrt{1-||B||^2}}=\frac{1-A\cdot B}{\sqrt{1-||A||^2}\sqrt{1-||B||^2}}.
$$
The equations give four symmetric solutions. One of them is:  $A(0.797 , 0)$, $B(0.496,  0.623)$. Then we compute the hyperbolic distance between all three points using the same formula and it will be $\sim$1.0905. This is the edge distance of a hyperbolic triangle that tiles the hyperbolic plane with seven triangles at each vertex. 

In a similar way one can compute edge length of an equilateral triangle with eight edges. Its edge is $\sim$1.5285 and each angle is $\pi/4 = 45^\circ$.

The above calculations demonstrate a fact that there are no similar triangles in the hyperbolic plane:{\it  The angles of a hyperbolic triangle determine its sides. }

Figure 4 shows only few layers of triangles around the origin. It appears triangles approach the boundary of the unit disc very fast but this process is infinite and there are infinitely many triangles in the tiling of a hyperbolic plane. There are four such layers in the crochet model in figure \ref{flower3}. 

\subsection{Circumference of a disk}

The ratio of a circumference of a disk to its diameter is constant and equal to $\pi$, remind that $C_{E^2}(r)=2\pi r$ in Euclidean plane. In hyperbolic plane the formula for the circumference of a disk with the same radius looks similar $C_{H^2}(r)=2 \pi \sinh{r}$ \cite{F}. The comparison of how fast $C_{H^2}(r)$ grows with respect to its diameter shows an exponential growth, see figure 5. This fact is also visible on the hyperbolic crochet models. We see the white and green boundaries on our models are curving in space and cannot be flatten.

\begin{figure}[h]
  \begin{center}
  \includegraphics[scale=0.7]{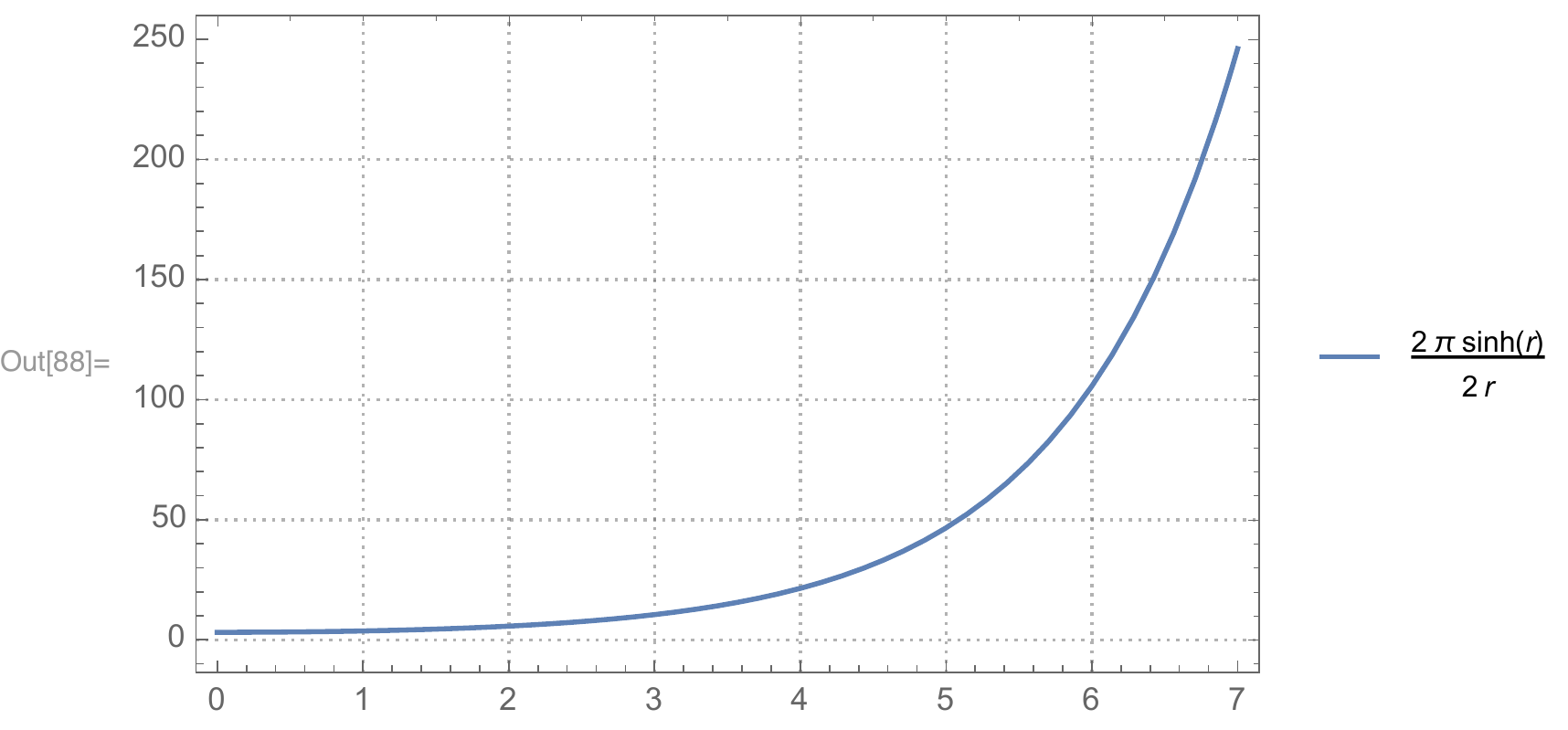}
  \end{center}
  \vspace{-0.3cm}
  \caption{Ratio of Circumference growth}
\end{figure}

\subsection{$\mathcal{C}^1$ embedding}

Our models have finite diameter and it is obvious that they occupy some finite size ball in $\mathbb R^3$.  We want to understand why they look so symmetrically "curvy" and the bigger the radius the more they look like curly balls. 

\medskip

Nash-Kuiper  $\mathcal{C}^1$-embedding theorem states that if a compact Riemannian $n$-manifold has a $\mathcal{C}^1$-embedding into $\mathbb R^{n+1}$ then it has an isometric embedding into $\mathbb R^{n+1}$ \cite{Kuiper,Nash}. As a consequence of this theorem it follows that a closed oriented Riemannian surface can be  $\mathcal{C}^1$-isometrically embedded into an arbitrarily small $\epsilon$-ball in Euclidean 3-space. And, there exist $\mathcal{C}^1$-isometric embeddings of the hyperbolic plane in $\mathbb R^3$. This says we can isometrically embed hyperbolic disk of finite radius into arbitrarily small ball $B_\epsilon$ in  $\mathbb R^3$. This is true only if the surface does not have any thickness and has zero volume. In practice, it can not be embedded into a ball of arbitrarily small radius. The discrete version of this theorem was proven by Burago-Zalgaller \cite{BZ} and implies that a 2-dimensional surface can be isometrically embedded into $\mathbb R^3$. Our crochet models with triangulation can be considered as piece-wise linear surfaces and can be used as a visual example of this theorem.

Next we want to understand why these models look curvy this particular way. The boundary can be thought as a rope of finite length and volume. The embedding ball in $\mathbb R^3$ has a finite volume and we fill it with a thick surface with some non-zero volume. As radius grows the circumference and the volume of a tube around the boundary grow exponentially but volume of a ball in $\mathbb R^3$ grows linearly. That is the reason why models of hyperbolic disks look less curvy and locally more flat (or smashed) but as radius becomes bigger they take more space of the ball, and after several steps volumes of the model surface and the ball get close. Then the hyperbolic disk turns into a ball and the we can not proceed further with this construction.
 
\section{Cylinder in a hyperbolic space}

How does a cylinder look like in a hyperbolic space? Note that hyperbolic cylinder is defined as a cylindrical surface in Euclidean space made of parallel lines passing through a hyperbola. There are different equivalent definitions of a cylinder in the Euclidean space, not all of them can be used in the hyperbolic space. One that works states that the {\it cylinder} in $\mathbb H^2$ is a surface of revolution obtained by a rotation of the line parallel to the axis of rotation.  Cutting the (infinite) cylinder along that geodesic line gives an (infinite) strip with two parallel lines. They have a unique common perpendicular which follows from the properties of hyperbolic geometry. It is called the {\it neck of a cylinder}.

\begin{figure}[h]
  \begin{center}
  \includegraphics[scale=0.38]{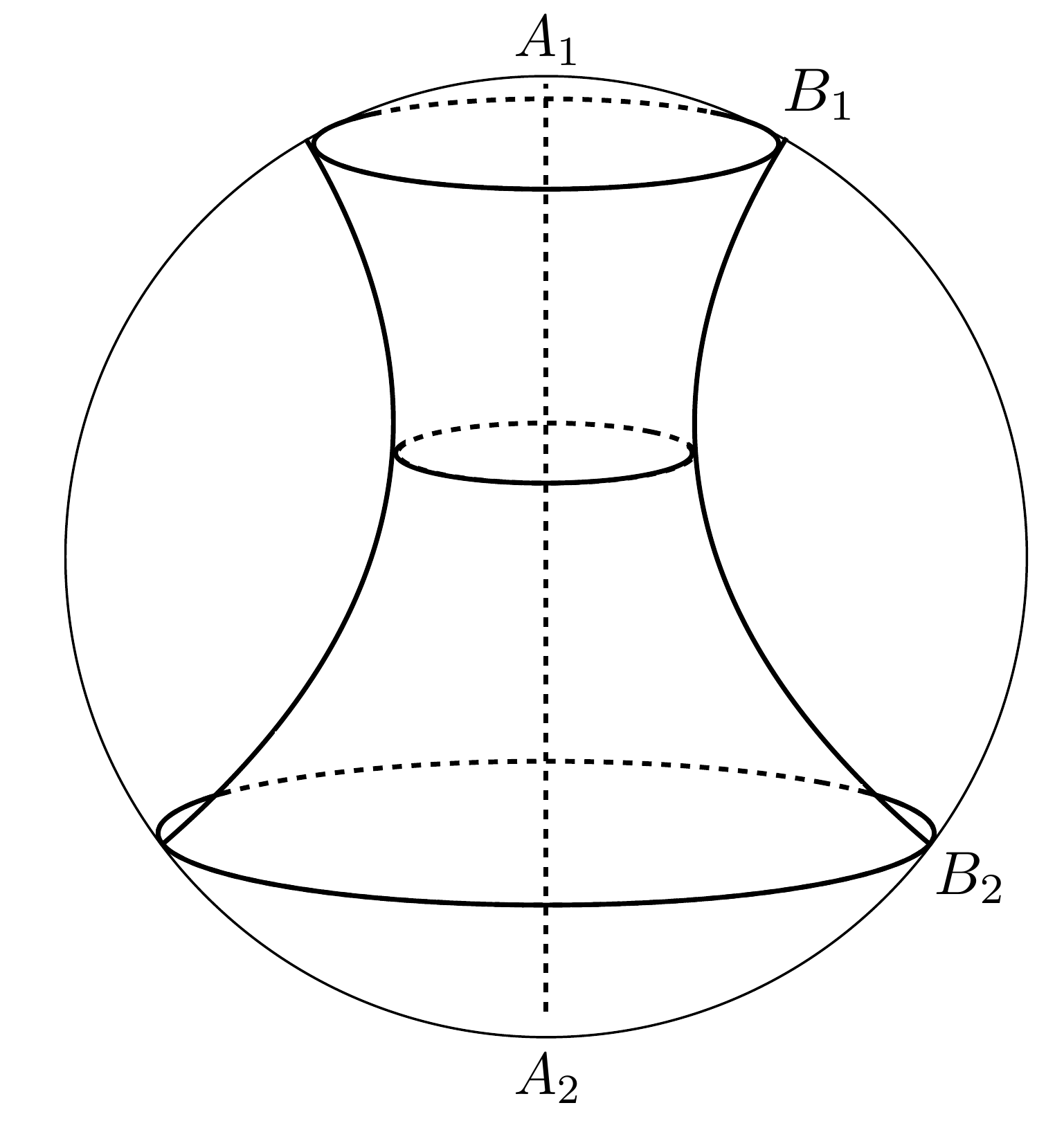}
  \vspace{-0.2cm}
  \caption{An infinite cylinder in the Poincare ball model.}
  \label{cylinder}
   \end{center}
\end{figure}

It is better to see a cylinder in the Poincare model of a hyperbolic space $\mathbb H^3$ which is a unit ball without its boundary sphere. Geodesics in this model are circular arcs perpendicular to the boundary of $\mathbb H^3$ and/or a chord passing through the origin.

In the figure 6 the axis of rotation $A_1 A_2$ is fixed and the geodesic $B_1 B_2$ is rotated. It forms a cylinder with two ideal circles on the boundary of $\mathbb H^3$. There is a unique plane perpendicular to both the axis $A_1A_2$ and the cylinder. Calculations in this model tell us that this plane intersects the cylinder along a geodesic which looks like a circle. This is the only geodesic on the infinite cylinder and it is the shortest simple closed curve on the cylinder.  

\vspace{-0.5cm}

\begin{figure}[h]
  \begin{center}
  \includegraphics[scale=0.2]{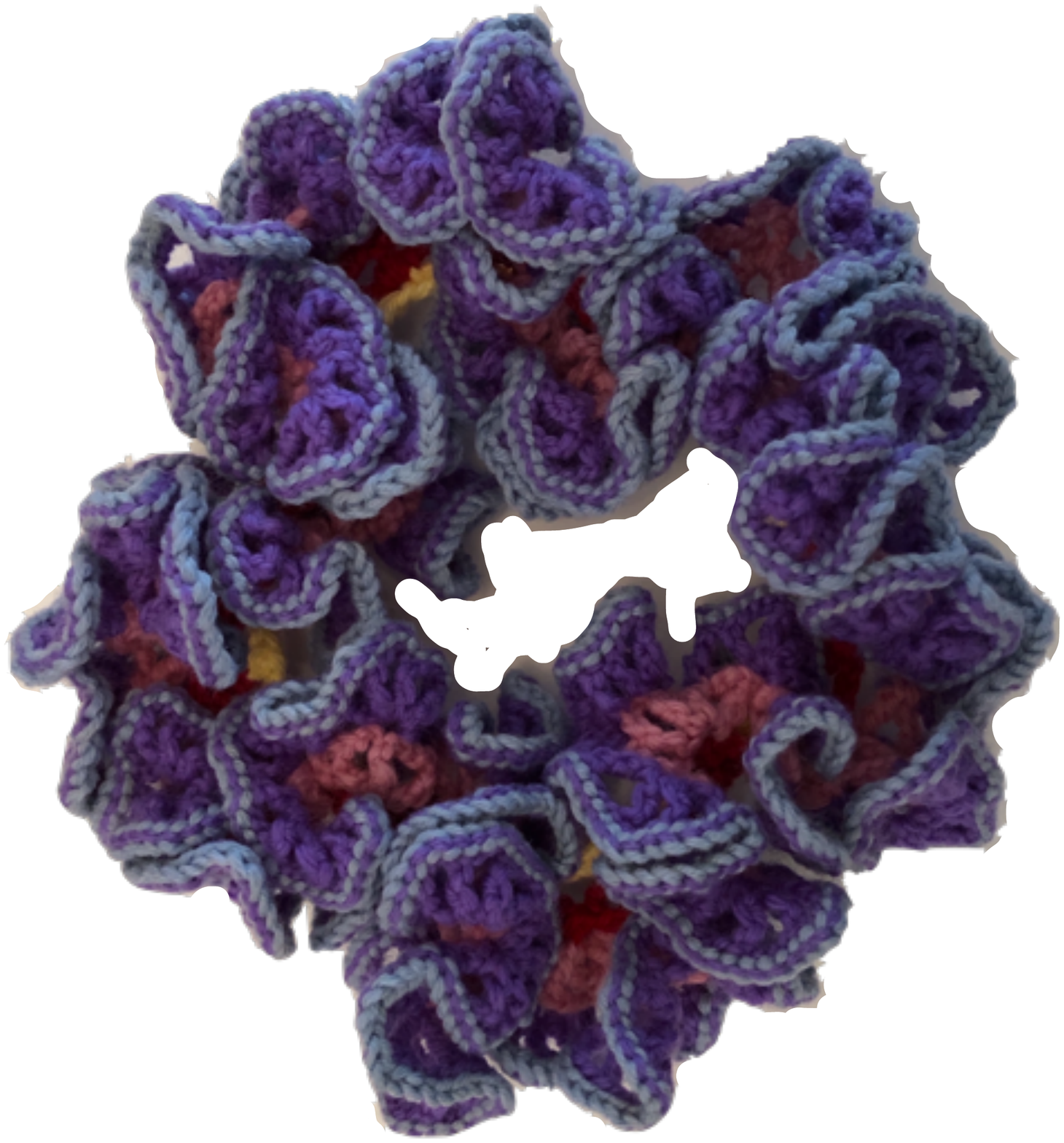}\quad
  \includegraphics[scale=0.2]{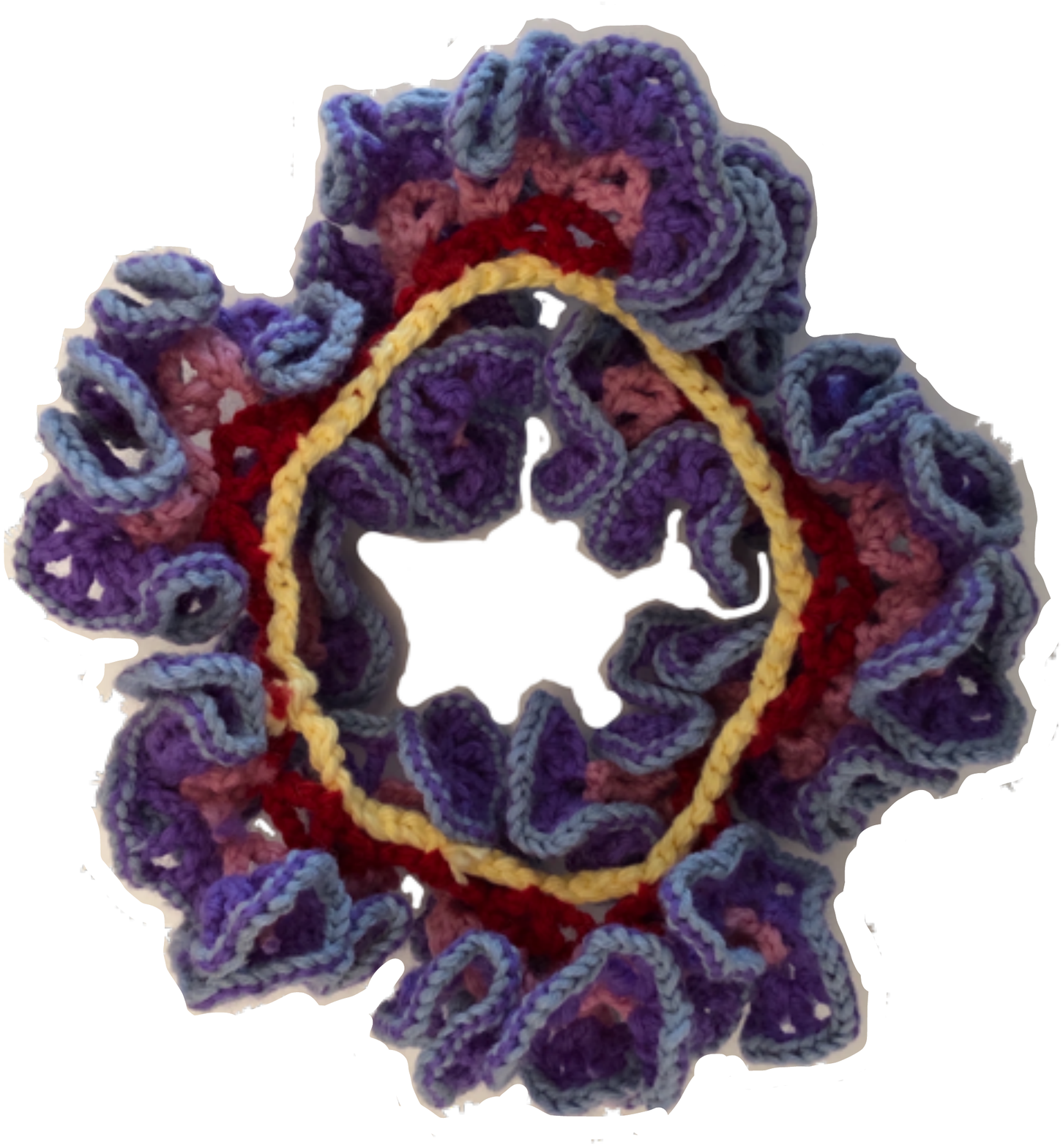}
  \vspace{-0.5cm}
  \caption{Part of an infinite cylinder bounded by its neck. Top and bottom views.}
   \end{center}
\end{figure}

A crochet model of a cylinder  in figure 7 is made out of eight equilateral triangles at each vertex.  It is bounded by its neck on one side (made with yellow yarn). The other boundary is only an approximated circle, it is formed by edges of equilateral triangles. It is contoured with blue yarn to stretch the edges. It also makes triangles closer to equilateral and the model looks overall nicer.

\section{Crochet and techniques}
The Gaussian curvature $K$  of a surface at a point is defined as the product of its two principal curvatures at the given point. Although the  curvature of the hyperbolic plane is a negative constant, the  crochet models are only its approximations as they have thickness along edges and empty at  the interior of triangles. The curvature of the models is discrete at vertices  and along edges of triangles. 

The way how the models look depend on yarn properties. More silky and stretchy yarn makes models that do not hold their shapes well (mushy). Yarn without any stretching like cotton is very hard to crochet tight and again it does not hold shape as we want. Rough yarn with a little bit of stretching allows to make rigid models (firm) with beautiful patterns of embedding. Wool or acrylic works well to maintain these properties. It has small stretching which allows to see local flatness of triangles and at the same time to hold the shape of an embedded surface in $\mathbb R^3$. 

The pattern for the seven triangles case is given on figure 9. For other models similar patterns were used.  
Crochet begins with a slip knot. 

\medskip

\begin{figure}[h]
\begin{center}
    \includegraphics[scale=0.1]{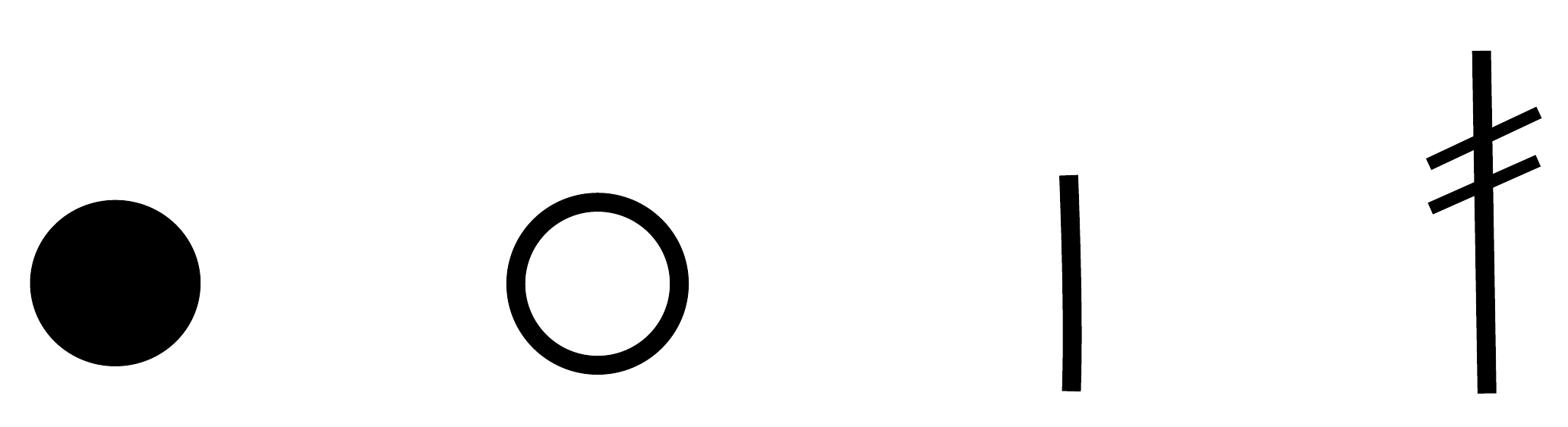}
        \vspace{-0.3cm}
  \caption{Crochet symbols (left to right): slip stitch, chain, half single and double crochet.}
\end{center}
\end{figure}

\medskip

\noindent {\bf Row 1:} 
\begin{itemize}
\item[Step 1:]  Crochet three chain stitches for rising, 
\item[Step 2:]   crochet three chain stitches, then double crochet stitch  in the first chain, 
\item[Step 3:] repeat previous step five more times, 
\item[Step 4:]  crochet three chains, 
\item[Step 5:] finish the row with the half basic stitch. 
\end{itemize}

\begin{figure}[h]
\begin{center}
    \vspace{-0.7cm}
        \includegraphics[scale=0.25]{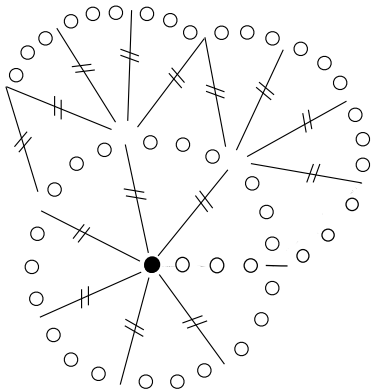}\\
        \vspace{-0.5cm}
  \caption{Crochet scheme for seven triangles at each vertex.}
\end{center}
\end{figure}

\newpage 

\noindent{\bf Row 2:} 
\begin{itemize}
\item[Step 6:] Three chain stitches for rising, 
\item[Step 7:] three chain stitches, double crochet stitch  into the first double chain from the previous row, 
\item[Step 8:] repeat the previous step three more times into the same double chain,   
\item[Step 9:] repeat steps 7$-$8 five more times. 
\item[Step 10:] three chains, 
\item[Step 11:]  finish the row with the half basic stitch. 
\end{itemize}

\medskip

\noindent{\bf Rows 3,4:} 

\noindent Repeat the pattern of Row 2. Yarn may be changed for another color when rows change.

\vspace{0.5cm}

\begin{figure}[h]
\begin{center}
    \includegraphics[scale=0.17]{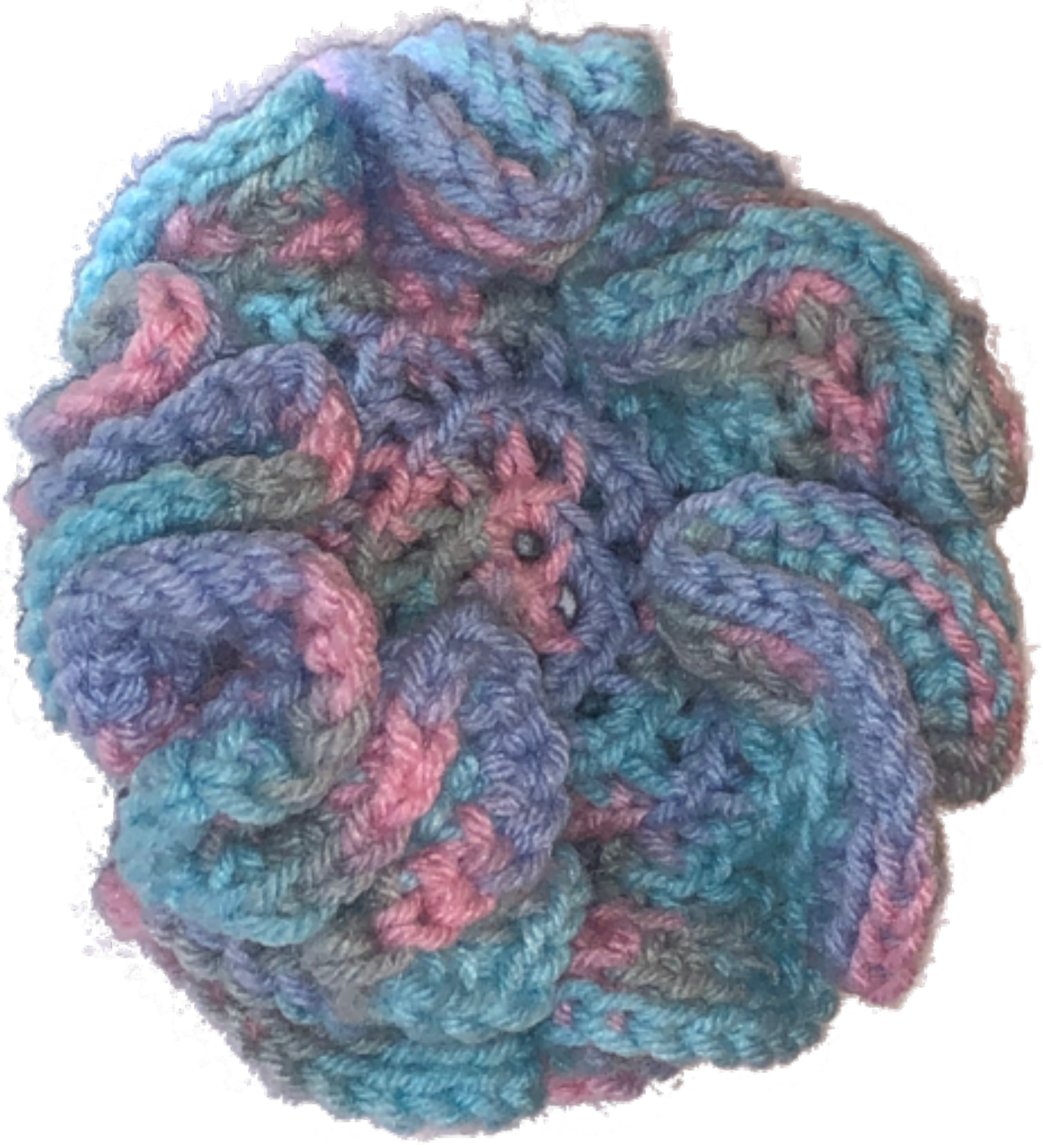}
    \includegraphics[scale=0.18]{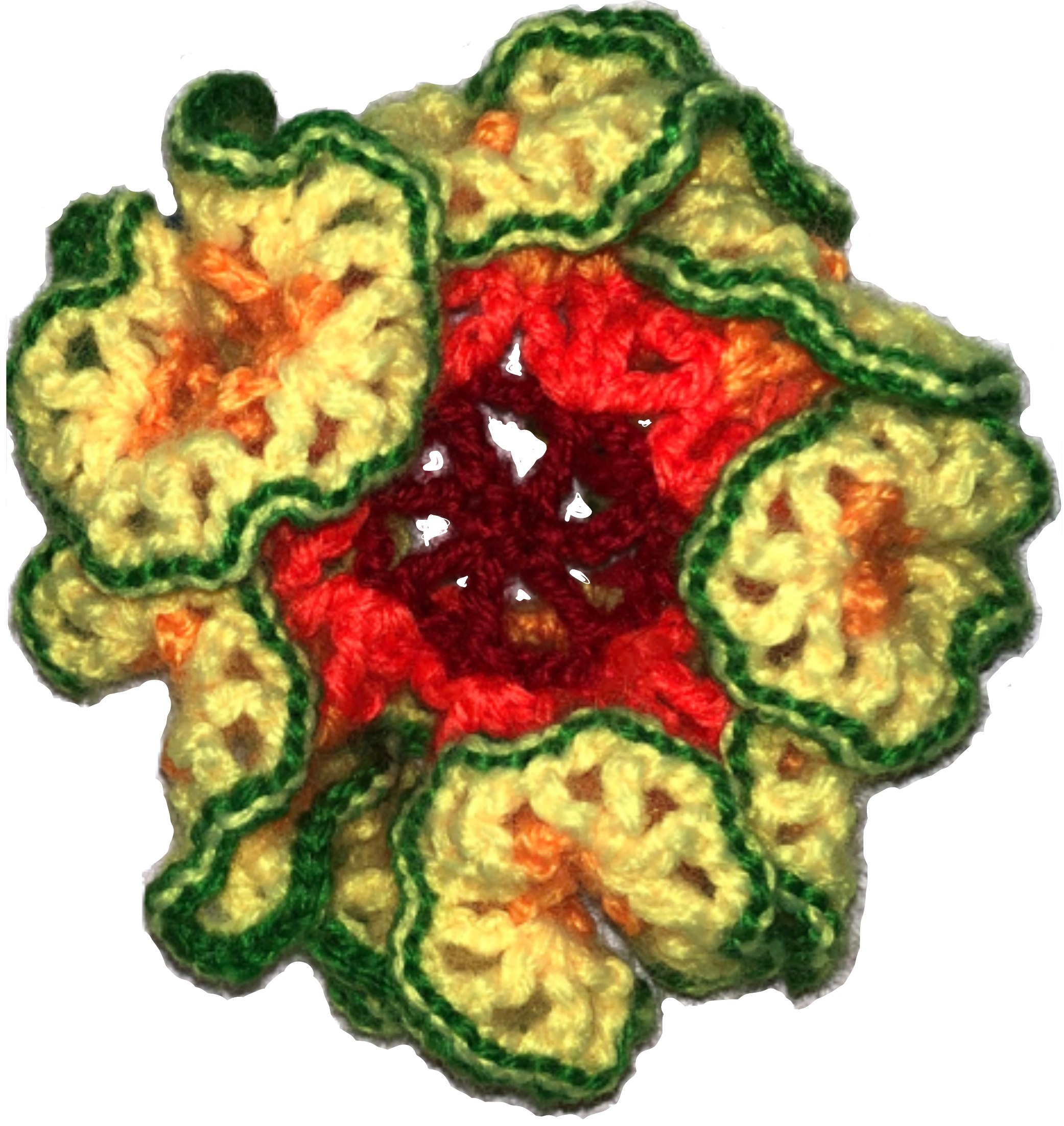}
        \includegraphics[scale=0.18]{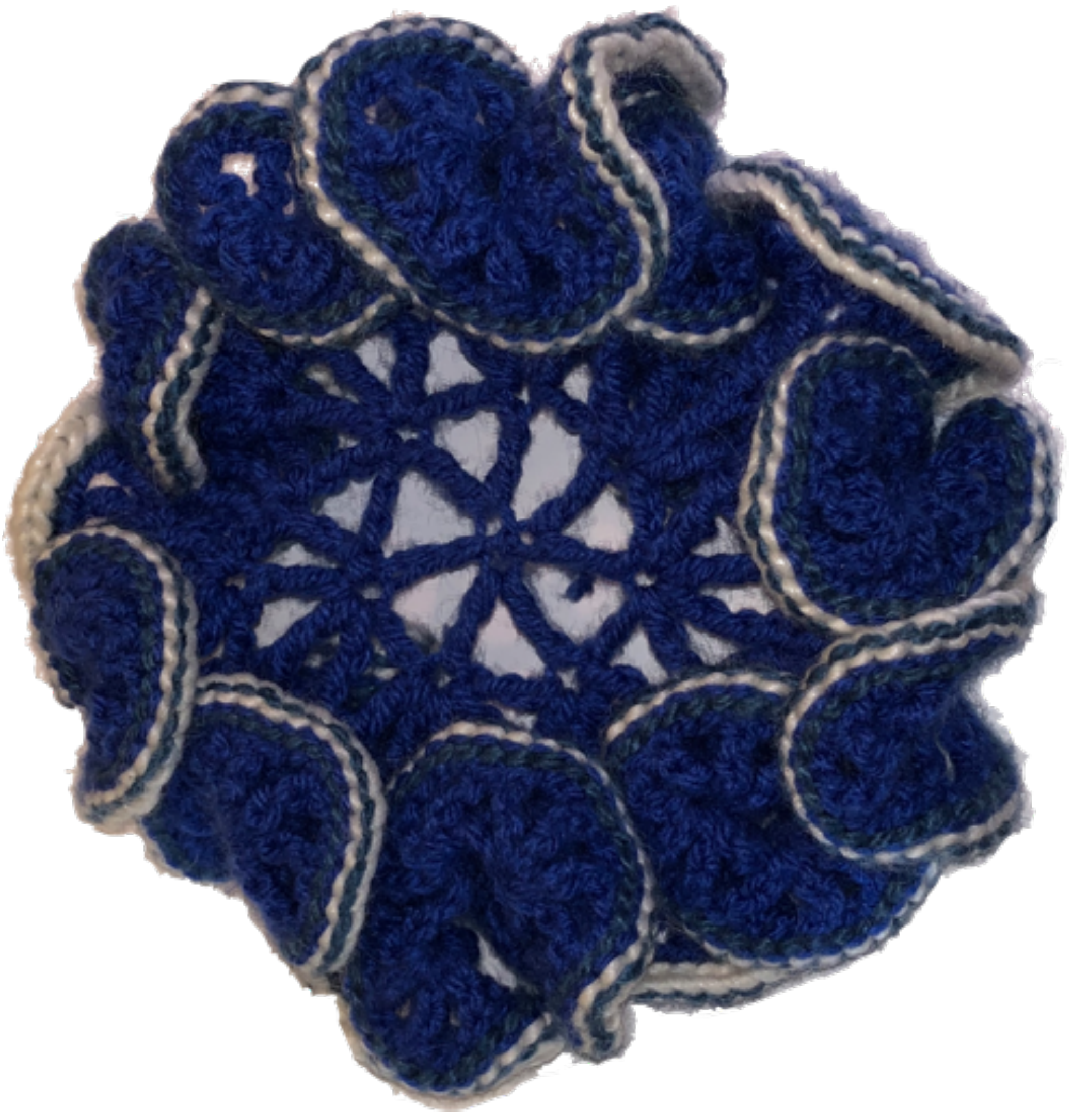}\\
        \vspace{-0.2cm}
  \caption{This is a closer look of the models. It is much more pleasant to appreciate these models in dimension 3  than from a picture of them.}
\end{center}
\end{figure}

\section{Artist's statement}

Maria Trnkova is a Krener Assistant Professor at University of California in Davis. She made crochet models of several hyperbolic surfaces to visualize them and explain certain features of hyperbolic geometry. In the process of making crochet models she used wool and acrylic yarn. She hopes her work would help others to learn hyperbolic geometry through visualization of these models. She was inspired by crochet models by Daina Taimina and the new models should be a useful complement for Taimina's models. 

Maria's personal story of discovering triangulated crochet models of hyperbolic disks started in the graduate class at Caltech. She was making preparations for the Riemannian Geometry class and wanted to demonstrate students how a hyperbolic disk looks like and visualize some of the geometric concepts. The process of crochet a solid hyperbolic disk is exponential in time and yarn consumption. She wanted to save time and use less yarn and as a trained crochet maker then experimented with different stitches and techniques. Very quickly the triangular pattern came up naturally.

\medskip

{\bf Acknowledgement:} The author is thankful to  Henry Segerman for encouragement in writing this note, to John Sullivan for asking a motivating question why these models look this particular way and to Joel Hass for an idea of depth coloring.

\end{document}